\pdfoutput=1

\documentclass[12pt]{article}
\usepackage[utf8]{inputenc}
\usepackage{amsmath}
\usepackage{amsthm}
\usepackage{amscd}
\usepackage{comment}
\usepackage[english]{babel}
\usepackage{color,graphicx,amssymb,tikz-cd,wrapfig}
\usepackage{amsmath,amsthm,amsfonts,amssymb,stmaryrd}
\usepackage{bm}
\usepackage[rightcaption]{sidecap}
\usepackage[margin=20truemm]{geometry}

\usepackage{enumerate}
\usepackage{pb-diagram}
\usepackage{graphicx} 
\usepackage[all]{xy}
\usepackage{appendix}
\title{Controlled $K$-theory and $K$-homology}
\author{Ryo Toyota}
\date{}

\begin{document}

\newtheorem{dfn}{Definition}[section]
\newtheorem{thm}{Theorem}[section]
\newtheorem{lem}{Lemma}[section]
\newtheorem{cor}{Corollary}[section]
\newtheorem{rmk}{Remark}[section]
\newtheorem*{pf}{Proof}
\newtheorem{ex}{Example}[section]
\newtheorem{prop}{Proposition}[section]
\newtheorem*{claim}{Claim}

\newcommand{\CC}{ {\mathbb C} }
\newcommand{\HH}{{\mathcal H}}
\newcommand{\FF}{{\mathbb F}}
\newcommand{\ZZ}{\mathbb{Z}}
\newcommand{\KK}{\mathcal{K}}
\newcommand{\diag}{\text{diag}}

\maketitle

\begin{abstract}
Motivated by the idea that our access to the spacetime is limited by the resolution of our measuring device, we give a new description of $K$-homology with a finite resolution. G. Yu introduced a $C^*$-algebra called the localization algebra $C^*_L(X)$ which consists of functions from $[1,\infty)$ to the Roe algebra $C^*(X)$ whose propagations converge to $0$ and he showed that for any finite dimensional simplicial complex $X$ endowed with the spherical metric, the $K$-theory of the localization algebra is isomorphic to the $K$-homology of $X$. We give a coarse graining version of this theorem using controlled $K$-theory (also known as quantitative $K$-theory). Namely, instead of considering families of operators whose propagations converge to $0$, we prove that for each dimension $n$, there exists a threshold $r_n>0$ such that the $K$-homology of $n$-dimensional finite simplicial complex $X$ is isomorphic to a certain group of equivalence classes of operators whose propagation is less than $r_n$. This picture also enables us to represent any element in the $K$-homology group $K_*(X)$ by a finite matrix for a finite simplicial complex $X$.
\end{abstract}

\tableofcontents

\section{Introduction}

In this paper, we introduce a quantitative picture of the $K$-homology group $K_*(X)$ for a finite simplicial complex $X$, which enables us to present every element in the $K$-homology group by a finite dimensional matrix following the idea of coarse graining in physics as discussed in \cite{connes2022tolerance}. In \cite{connes2022tolerance}, based on the idea that we can only determine the underlying metric space up to a finite resolution, A. Connes and W. D. Suijlekom proposed a noncommutative geometric framework to encode our limited resolution to physical measurements using a tolerance relation (symmetric and reflexive but not transitive relation). We approach $K$-homology based on a tolerance relation $\mathcal{R}=\{(x,y):d(x,y)< r\}\subset X\times X$ for a fixed $r>0$.

The $K$-homology group has already several well-known descriptions. Its first picture was given by Kasparov in \cite{kasparov1975topological}. He introduced the notion of the Fredholm modules and defined $K$-homology to be a certain group of equivalence classes of them, which is a topological invariant. The idea of the Fredholm module originated with a functional analytic abstraction of elliptic operators by Atiyah \cite{atiyah1970global}.

Another picture was given by G. Yu in \cite{yu1997localization}. It is known that the $K$-theory of Roe algebras $C^*(X)$ depends only on the large scale structure of $X$, so it is not isomorphic to $K$-homology in general. But Yu defined a $C^*$-algebra $C^*_L(X)$ called the localization algebra which consists of functions $f:[1,\infty) \rightarrow C^*(X)$ such that the propagation of an operator $f(t)$ converges to $0$. The idea behind it is that if we restrict our attention to operators whose propagations are small, then we can recover local structures of $X$ which we lost by forming the Roe algebra. Actually, he showed that the $K$-theory of localization algebra is isomorphic to the $K$-homology for finite dimensional simplicial complexes. In \cite{qiao2010localization}, Y. Qiao and J. Roe showed that this holds for any locally compact metric space if the module is very ample. 
The goal of this paper is to realize $K$-homology using the controlled $K$-theory groups $K^{\epsilon,r}_*(C^*(X))$. Instead of considering families of projections or unitaries whose propagations converge to 0, we only consider projections and unitaries whose propagations are smaller than a fixed threshold $r$. This is the coarse graining picture of $K$-homology in the same spirit as \cite{connes2022tolerance}.
The recent progress to approximate $KK$-theory of $C^*$-algebras by controlled $K$-theory can be found \cite{willett2020controlled}, but this paper focuses more on elementary examples ($K$-homology of finite simplicial complex) to give a description by finite matrices.

For some technical flexibilities, we introduce the notions of quasi-projection and quasi-unitary.
Throughout this paper, an operator $p$ over $C^*(X)$ is called an $(\varepsilon,r)$-quasi-projection if the propagation of $p$ is less than $r$ and $p$ satisfies $\|p^2-p\|\leq \varepsilon$ and $p=p^*$. An operator $u$ over $C^*(X)$ is called an $(\varepsilon,r)$-quasi-unitary if the propagation of $u$ is less than $r$ and $u$ satisfies $\|u^*u-1\|\leq \varepsilon$ and $\|uu^*-1\|\leq \varepsilon$. The controlled $K$-group $K^{\varepsilon,r}_*(C^*(X))$ is defined to be the group of equivalence classes of $(\varepsilon,r)$-quasi-unitaries over $C^*(X)$, where the equivalence relation is the homotopy relation throughout $(3\varepsilon,2r)$-quasi-unitaries for $*=1$. The more detailed definition including the case for $*=0$ can be found in Definition \ref{quantitative k-theory}. For $0< \varepsilon < \varepsilon'<\frac{1}{4}$ and $0<r<r'$, we have a canonical map, which is called the forgetful map
\begin{align*}
    \iota^{(\varepsilon,r),(\varepsilon',r')}:K^{\varepsilon,r}_*(C^*(X)) \rightarrow K^{\varepsilon',r'}_*(C^*(X)).
\end{align*}
The image of this map $\iota^{(\varepsilon,r),(\varepsilon',r')}(K^{\varepsilon,r}_1(C^*(X)))$ can be regarded as a group generated by $(\varepsilon,r)$-quasi-unitaries with a relaxed equivalence relation, homotopy throughout $(3\varepsilon',2r')$. We denote this group by $K_*^{(\varepsilon,r),(\epsilon',r')} (C^*(X))$.
 Then we can state the main theorem.

\begin{thm}\label{main theorem}
  \rm{For each dimension $n$, there exist a constant $\lambda_n$, a function $h_n:(0,\frac{1}{4\lambda})\rightarrow [1,\infty)$, constants $r_n$ and $\varepsilon_n$ depending only on $n$ such that for any $n$-dimensional finite simplicial complex $X$, we have
    \begin{align}\label{main}
        K_*(X) \cong K_*^{(\varepsilon,r),(\lambda_n\epsilon,h_n(\varepsilon)r)} (C^*(X))   
    \end{align}
    for all $(\varepsilon,r)$ with $0< \varepsilon< \varepsilon_n$ and $0<r<r_n$.
    }
\end{thm}

The main technical tool of the proof of Theorem \ref{main theorem} is an asymptotically exact Mayer-Vietoris sequence \cite{oyono2019quantitative} to extract information from symplicial complexes whose dimension is smaller than that of $X$. Note that the main idea of the proof appeared \cite{yu2010characterization}, but our main theorem should be the precise statement of Theorem 4.7 of \cite{yu2010characterization}.

If $X$ is compact, then its Roe algebra is the set of all compact operators, which is the inductive limit of matrix algebras. So each $K$-homology element of a finite simplicial complex can be expressed by a finite matrix via the isomorphism \eqref{main}, as explained in Remark \ref{discretization} with a geometric picture of discretization.

\section{Roe algebras and controlled $K$-theory}

In this section, we define some basic concepts such as support of operators, propagation, Roe algebra, filtered $C^*$-algebra, and controlled $K$-theory. 

For any Hilbert space $H$, we denote the set of all compact operators by $\KK(H)$.

\begin{dfn}\label{ample repn}
\rm{    Let $X$ be a locally compact metric space and $H_X$ a separable Hilbert space. We say $H_X$ is an ample $X$-module if there is a non-degenerate $*$-homomorphism
    \begin{align*}
        \rho:C_0(X) \rightarrow B(H_X)
    \end{align*}
    such that $\rho(f)$ is not a compact operator for any $f \in C_0(X) \setminus \{0\}$.
    When it is clear from the context, we will write $f$ in place of $\rho(f)$.
    }
\end{dfn}

\begin{dfn}\label{support}
\rm{
    Let $H_X$ be an $X$-module and $H_Y$ be an $Y$-module. For $T \in B(H_X,H_Y)$, we define a subset $\text{supp}(T) \subset Y\times X$ to be the complement of the set of all $(y,x)\in Y\times X$ such that there exist two functions $f\in C_0(X)$ and $g \in C_0(Y)$ such that
    \begin{align*}
        gTf = 0, \text{ } f(x) \neq 0, \text{ } g(y) \neq 0.
    \end{align*}
   Also we define the propagation of $T$ by
   \begin{align*}
       \text{prop}(T):=\sup\{d(x,y);(y,x) \in \text{supp}(T)\} \in [0,\infty].
   \end{align*}
}
\end{dfn}

 We now define a $C^*$-algebras $C^*(X)$, which will be used throughout this paper. 
 
\begin{dfn}\label{Roe algebra}
\rm{Let $H_X$ be an ample non-degenerate $X$-module. We define the following algebras:
\begin{align*}
    C^*_{\text{alg}}(H_X):=\{T\in B(H_X);Tf \in \KK(H_X) \text{ for any }f \in C_0(X) \text{ and } \text{prop}(T)<\infty\} 
\end{align*}
\begin{align*}
    C^*(H_X):=\overline{C^*_{\text{alg}}(H_X)},
\end{align*}
  where the closure is taken with respect to the operator norm in $B(H_X)$.
    
When the module $H_X$ is clear from the context, we will write $C^*(X)$ instead of $C^*(H_X)$.
}    
\end{dfn}

\begin{dfn}\label{filter}
\rm{    A filtered $C^*$-algebra $A$ is a $C^*$-algebra equipped with a family $(A_r)_{r>0}$ of closed linear subspaces parametrized by positive numbers $r$ such that:

    \begin{enumerate}[1)]
        \item $A_r \subset A_{r'}$ if $r \leq r'$

        \item $A_r^*=A_r$

        \item $A_r\cdot A_{r'} \subset A_{r+r'}$

        \item the subalgebra $\cup_{r>0} A_r$ is dense in $A$.
    \end{enumerate}
    
If $A$ is unital we also require that the identity $1$ is an element of $A_r$ for every positive number $r$.
}
\end{dfn}

\begin{rmk}\label{closure}
\rm{
    A Roe algebra $C^*(X)$ is filtered by propagation of operators, namely $C^*(X)_r$ is the closure of the set of operators whose propagation is smaller than $r$.
}
\end{rmk}

We can define the controlled $K$-theory for filtered $C^*$-algebras, where ``control" refers to control of propagations. 

\begin{dfn}\label{quasi-projectio}
\rm{    Let $A=(A_r)_{r>0}$ be a unital filtered $C^*$-algebra and $\varepsilon$ and $r$ be positive numbers. We define the sets of $(\varepsilon,r)$-quasi-projections and $(\varepsilon,r)$-quasi-unitaries in $A$ respectively by
    \begin{align*}
        P^{\varepsilon,r}(A)&:=\{p \in A_r: p^*=p,\|p^2-p\|<\varepsilon\}\\
        U^{\varepsilon,r}(A)&:=\{u \in A_r:\|u^*u-1\|<\varepsilon,\|uu^*-1\|<\varepsilon\}.
    \end{align*}
    For any positive integer $n$, we set $P^{\varepsilon,r}_n(A)=P^{\varepsilon,r}(M_n(A))$ and $U^{\varepsilon,r}_n(A)=U^{\varepsilon,r}(M_n(A))$. Then there are canonical inclusions
    \begin{align*}
        P^{\varepsilon,r}_n(A) \hookrightarrow P^{\varepsilon,r}_{n+1}(A);p \mapsto \begin{pmatrix}
            p & 0 \\
            0 & 0
        \end{pmatrix}
    \end{align*}
    and
    \begin{align*}
        U^{\varepsilon,r}_n(A) \hookrightarrow U^{\varepsilon,r}_{n+1}(A);u \mapsto \begin{pmatrix}
            u & 0 \\
            0 & 1
        \end{pmatrix}.
    \end{align*}
With respect to these inclusions, their unions are denoted by
\begin{align*}
    P^{\varepsilon,r}_{\infty}=\cup_{n=1}^{\infty} P^{\varepsilon,r}_n
\end{align*}
and
\begin{align*}
    U^{\varepsilon,r}_{\infty}=\cup_{n=1}^{\infty} U^{\varepsilon,r}_n.
\end{align*}
}
\end{dfn}

\begin{dfn}\label{quantitative k-theory}
\rm{    Let $A=(A_r)_r$ be a filtered unital $C^*$-algebra, $r$ and $\varepsilon$ be positive numbers with $\varepsilon<1/4$. We define the following equivalence relations on $P_{\infty}^{\epsilon,r}(A)\times \mathbb{N}$ and $U^{\epsilon,r}_{\infty}(A)$ .
\begin{enumerate}[1)]
    \item  For $p,q\in P^{\varepsilon,r}_{\infty}(A)$ and $\ell,\ell'\in \mathbb{N}$, $(p,\ell)\sim_{\varepsilon,r}(q,\ell')$ if there exist a positive integer $k$ such that $\text{diag}(p,I_{k+\ell'})$ and $\diag(q,I_{k+\ell'})$ are homotopic in $P^{\varepsilon,r}_{\infty}(A)$.

    \item For $u,v\in U^{\varepsilon,r}_{\infty}(A)$, $u\sim_{3\varepsilon,2r}v$ if $u$ and $v$ are homotopic in $U_{\infty}^{3\varepsilon,2r}(A)$.

   By these equivalence relations, the controlled $K$-groups are defined by
    \begin{align*}
        K_0^{\varepsilon,r}(A)&:=(P^{\varepsilon,r}_{\infty}(A) \times \mathbb{N}) / \sim_{\varepsilon,r}\\
        K_1^{\varepsilon,r}(A)&:=U^{\varepsilon,r}_{\infty}(A) / \sim_{3\varepsilon,2r}.
    \end{align*} 
    If $A$ is non unital, by letting $\Tilde{A}=(A_r+\mathbb{C})_r$ be the unitization of $A$, define
    \begin{align*}
        K_0^{\varepsilon,r}(A)&:=\{[p,\ell]_{\varepsilon,r} \in P^{\varepsilon,r}_{\infty}(\Tilde{A}) \times \mathbb{N}/\sim_{\varepsilon,r}:\text{ }\chi(\rho_A(p))=\ell\}\\
        K_1^{\varepsilon,r}(A)&:=U^{\varepsilon,r}_{\infty}(\Tilde{A}) / \sim_{3\varepsilon,2r},
    \end{align*}
    where $\rho_A:\Tilde{A}\twoheadrightarrow \mathbb{C}$ is the projection onto the scalar and $\chi$ is the characteristic function of the interval $[\frac{1}{2},\infty)$. Also we define the additive structures by
    \begin{align*}
        [p,\ell]_{\varepsilon,r}+[q,\ell']_{\varepsilon,r}&=[\diag(p,q),\ell+\ell']_{\varepsilon,r}\\
        [u]_{3\varepsilon,2r}+[v]_{3\varepsilon,2r}&=[\diag(u,v)]_{3\varepsilon,2r}.
    \end{align*}
\end{enumerate}
}
\end{dfn}

\begin{rmk}
\rm{
    With these notations, $K_*^{\varepsilon,r}(A)$ is an abelian group.
    The reason to consider $(3\varepsilon,2r)$-homotopy is to make $K_1$ an abelian group (Lemma 1.15 and Remark 1.17 of \cite{oyono2015quantitative}.)
    }
\end{rmk}

\begin{rmk}\label{almost same projections}
\rm{
    Let $A=(A_r)_r$ be a filtered unital $C^*$-algebra. For $0<\delta<\frac{1}{4}$, $0<\varepsilon<\frac{1}{4}$ and $r>0$, assume $p$ is an $(\varepsilon,r)$-quasi-projections. If $p'$ is an operator whose propagation is smaller than $r$ and $\|p-p'\|<\delta$, then $p'$ is a $(\varepsilon+5\delta,r)$-quasi-projection and the homotopy which connect $p$ and $p'$ linearly is a homotopy through $(\varepsilon+5\delta,r)$-quasi-projection. This is because we have
    \begin{align*}
        \|p'^2-p'\|&\leq \|p'^2-p'p\|+\|p'p-p^2\|+\|p^2-p\|+\|p-p'\|\\
        &\leq 2\delta+2\delta +\varepsilon +\delta=5\delta+\varepsilon
    \end{align*}
    and for any $t\in [0,1]$ by applying the above formula to $p_t=(tp+(1-t)p')$, 
    \begin{align*}  
        &\|p_t^2-p_t\|\leq 5\|p-p_t\|+\varepsilon\leq 5\delta+\varepsilon.
    \end{align*}
    The corresponding statement is also true for quasi-unitaries.
    }
\end{rmk}

\begin{dfn}
\rm{    We have the following maps from the controlled $K$-theory to $K$-theory. Assume $A$ is a unital filtered $C^*$-algebra.
    
    Let $\chi=\chi_{[\frac{1}{2},\infty)}$ be a characteristic function of $[\frac{1}{2},\infty)$. If $p\in M_n(A)$ is selfadjoint and satisfy $\|p^2-p\|<1/4$, then $\frac{1}{2} \notin \sigma(p)$. By continuous functional calculus, we have a map
    \begin{align*}
        \kappa : K_0^{\varepsilon,r}(A) \rightarrow K_0(A); [p] \mapsto [\chi(p)].
    \end{align*}

    If $u \in M_n(A)$ satisfies $\|uu^*-1\|<\frac{1}{4}$ and $\|u^*u-1\|<\frac{1}{4}$, then $u$ is invertible. We define
    \begin{align*}
        \kappa : K_1^{\varepsilon,r}(A) \rightarrow K_1(A);[u] \mapsto [u/(u^*u)^{\frac{1}{2}}].
    \end{align*}
    These two maps are called the comparison maps just after Definition 4.1 and Definition 4.4 of \cite{guentner2016dynamical}.
    }
\end{dfn}

\section{Overview of quantitative objects}
In this section, we introduce a general framework to deal with controlled $K$-theory as a family of groups $(K^{\varepsilon,r}_*(C^*(X)))_{\varepsilon,r}$ systematically. Everything up to Definition \ref{composition} is contained in the section 1 of \cite{oyono2019quantitative}. In Definition \ref{equivalence}, we introduce a notion of equivalence between two quantitative objects. By using this concept, for an $(n-1)$-dimensional finite simplicial complex $Z$ and an $n$-dimensional finite simplicial complex $\Tilde{Z}$, we can replace $K^{\varepsilon,r}_*(C^*(\Tilde{Z}))$ by $K^{\varepsilon,r}_*(C^*(Z))$ in a long exact sequence to use induction hypothesis on the dimension.

In the subsequent definitions we define quantitative objects and morphisms between them abstractly. Main examples of a quantitative objects are controlled $K$-theory of Roe algebras.
\begin{dfn}
\rm{      A quantitative object is a family $\mathcal{O}=(O^{\varepsilon,r})_{0<\varepsilon<1/4,r>0}$ of abelian groups, together with group homomorphisms
    \begin{align*}
        \iota ^{\varepsilon,\varepsilon',r,r'}:O^{\varepsilon,r} \longrightarrow O^{\varepsilon',r'} 
    \end{align*}
    for $0<\varepsilon \leq \varepsilon' <1/4$ and $0<r\leq r'$ such that
 \begin{enumerate}[1)]

     \item $\iota ^{\varepsilon,\varepsilon,r,r}=Id_{O^{\varepsilon,r}}$ for any $0<\varepsilon<1/4$ and $r>0$; 

     \item  $\iota ^{\varepsilon',\varepsilon'',r',r''}\circ \iota ^{\varepsilon,\varepsilon',r,r'}=\iota ^{\varepsilon,\varepsilon'',r,r''}$
 ,for any $0<\varepsilon \leq \varepsilon' \leq \varepsilon''<1/4$ and $0<r\leq r' \leq r''$.
\end{enumerate}

These maps $\iota ^{\varepsilon, \varepsilon', r, r'}$ are called structure maps.  }
    \end{dfn}

\begin{rmk}
\rm{
We apply abstract terminologies of quantitative object to $O^{\varepsilon,r}=K^{\varepsilon,r}_*(C^*(X))$ for locally compact metric space $X$. In this case, the structure maps $\iota^{\varepsilon,\varepsilon',r,r'}$ are the forgetful maps which assign an equivalence class $[p]$ to the equivalence class $[p]$ represented by the same operator $p$. 
}
\end{rmk}

\begin{dfn}
\rm{    A control pair is a pair $(\lambda,h)$ satisfying;

\begin{enumerate}[1)]
    \item $\lambda$ is a positive number with $\lambda>1$;

    \item $h:(0, \frac{1}{4\lambda}) \rightarrow (1,+\infty)$ is a map such that there exists a non-increasing map $g:(0, \frac{1}{4\lambda}) \rightarrow (1, +\infty) $, with $h \leq g$.
\end{enumerate}

For two control pair $(\lambda,h)$ and $(\lambda',h')$ we define its composition $(\lambda,h)*(\lambda',h')=(\lambda \lambda',h*h')$ by
\begin{align*}
    h*h':(0,\frac{1}{4\lambda\lambda'}) \rightarrow (1,+\infty ); \varepsilon \mapsto h(\lambda'\varepsilon)h'(\varepsilon)
\end{align*}

For a control pair $(\lambda,h)$ and a quantitative object $\mathcal{O}=(O^{\varepsilon,r})_{0<\varepsilon<1/4,r>0}$, we sometimes write its structure map as $\iota^{(\varepsilon,r),(\lambda,h)(\varepsilon,r)}$ in the sense of $\iota^{\varepsilon,\lambda \varepsilon, r, h(\varepsilon)r}$.
 }   
\end{dfn}

\begin{dfn}
 \rm{   Let $(\lambda,h)$ be a control pair, let $\mathcal{O}=(O^{\varepsilon,s})_{0<\varepsilon <1/4,s>0}$ and $\mathcal{O'}=(O'^{\varepsilon,s})_{0<\varepsilon <1/4,s>0}$ be quantitative objects and let $r$ be a positive number. A $(\lambda,h)$-controlled morphism of order $r$ 
    \begin{align*}
        \mathcal{F}:\mathcal{O} \rightarrow \mathcal{O'}
    \end{align*}
    is a family $\mathcal{F}=(F^{\varepsilon,s})_{0<\varepsilon<\frac{1}{4\lambda},0<s<\frac{r}{h(\varepsilon)}}$ of group homomorphisms
    \begin{align*}
        F^{\varepsilon,s}:O^{\varepsilon,s} \rightarrow O'^{\lambda \varepsilon,h(\varepsilon)s}
    \end{align*} which are compatible with structure maps.
    }
\end{dfn}

The next definition is about the exactness of a sequence of morphisms.
\begin{dfn}\label{composition}
\rm{    
Let $(\lambda,h)$ be a control pair and let $\mathcal{O}$, $\mathcal{O}'$ and $\mathcal{O}''$ be quantitative objects. Let
    \begin{align*}
        \mathcal{F}:\mathcal{O} \rightarrow \mathcal{O}'
    \end{align*} be a $(\alpha_{\mathcal{F}},k_\mathcal{F})$-controlled morphism and let
     \begin{align*}
        \mathcal{G}:\mathcal{O}' \rightarrow \mathcal{O}''
    \end{align*} be a $(\alpha_{\mathcal{G}},k_\mathcal{G})$-controlled morphism. Then the composition
    \[
\begin{CD}
    \mathcal{O} @>{\mathcal{F}}>> \mathcal{O}' @>{\mathcal{G}}>> \mathcal{O}''
\end{CD}
    \]
    is said to be $(\lambda,h)$-exact at $\mathcal{O}'$ of degree $r$ if $\mathcal{G}\circ \mathcal{F}=0$ and for any $0<\varepsilon<\frac{1}{4\max\{\lambda\alpha_{\mathcal{F}},\alpha_{\mathcal{G}}\}}$, any $0<s<\frac{1}{k_{\mathcal{F}}(h(\lambda\varepsilon))r}$ and any $y \in O'^{ \varepsilon, s}$ such that $G^{\varepsilon,s}(y)=0 $, there exists an element $x \in O^{\lambda \varepsilon, h(\varepsilon)s}$ such that
    \begin{align*}
        F^{\lambda\varepsilon,h(\lambda\varepsilon)s}(x)=\iota^{\varepsilon,\alpha_{\mathcal{F}}\lambda\varepsilon,s,k_{\mathcal{F}}(h(\lambda\varepsilon))s}(y),
    \end{align*}
    as the following diagram:
    \[
\begin{CD}
    @.y \in \mathcal{O}'^{\varepsilon,s} @>{G^{\varepsilon,s}}>> \mathcal{O}''^{\alpha_{\mathcal{G}}\epsilon,k_\mathcal{G}(\epsilon)s} \ni G^{\varepsilon,s}(y)=0\\
    @.  @V{\iota}VV @. \\
    x\in \mathcal{O}^{\lambda \varepsilon,h(\varepsilon)s} @>{F^{\lambda \varepsilon,h(\varepsilon)s}}>> \mathcal{O}'^{\alpha_{\mathcal{F}}\lambda\varepsilon,k_{\mathcal{F}}(h(\lambda\varepsilon))s }.@.\\  
\end{CD}
\]
    }
\end{dfn}

Next, we define the equivalence between quantitative objects.
\begin{dfn}\label{equivalence}
\rm{    
Let $(\lambda,h)$ be a control pair and $\mathcal{O}=(\mathcal{O}^{\varepsilon,r})_{0<\varepsilon<\frac{1}{4}}$ and $\mathcal{O}'=(\mathcal{O}'^{\varepsilon,r})_{0<\varepsilon<\frac{1}{4}}$ be quantitative objects. A $(\lambda,h)$-controlled morphism of degree $r$
    \begin{align*}
        \mathcal{F}:\mathcal{O} \longrightarrow \mathcal{O}'
    \end{align*} 
is said to be a $(\lambda,h)$-controlled equivalence morphism between $\mathcal{O}$ and $\mathcal{O}'$ if there exists a $(\lambda,h)$-controlled morphism
    \begin{align*}
        \mathcal{G}:\mathcal{O}' \longrightarrow \mathcal{O}
    \end{align*}
    such that for all $s>0$ with $h(\lambda \varepsilon)h(\varepsilon)s<r$, we have
    \begin{align*}
         G^{\lambda \varepsilon , h(\varepsilon)s}\circ F^{\varepsilon, s} &= \iota ^{\varepsilon, \lambda^2\varepsilon, s,h(\lambda \varepsilon)h(\varepsilon)s}\\
         F^{\lambda \varepsilon , h(\varepsilon)s}\circ G^{\varepsilon, s} &= \iota ^{\varepsilon, \lambda^2\varepsilon, s,h(\lambda \varepsilon)h(\varepsilon)s}.
    \end{align*}
    If such a control pair $(\lambda,h)$ and $(\lambda,h)$-controlled morphisms $\mathcal{F}$ and $\mathcal{G}$ exist, then we say that the two quantitative objects $\mathcal{O}$ and $\mathcal{O}'$ are asymptotically equivalent.
  }  
\end{dfn}

\begin{rmk}\label{closed or non-closed}
\rm{
    In some literature, each subspace $A_r$ is not required to be closed in Definition \ref{filter} and we can define the controlled $K$-theory in the same way without assuming it (for example \cite{oyono2015quantitative}).  Here we remark that this difference is not important $K$-theoretically. Assume $A=(A_r)$ be a unital filtered $C^*$-algebra with possibly ``non-closed" subspaces $A_r$ and we encode a different filtration on the same $C^*$-algebra by $A'=(A'_r)$ with $A'_r=\overline{A_r}$. Then one can show that two quantitative objects $(K^{\varepsilon,r}_*(A))_{\varepsilon,r}$ and $(K^{\varepsilon,r}_*(A'))_{\varepsilon,r}$ are asymptotically equivalent. We have a natural map
    \begin{align*}
        \iota^{\varepsilon,r}: K^{\varepsilon,r}_*(A)\rightarrow K^{\varepsilon,r}_*(A'); \text{ }[p,n]\mapsto [p,n].
    \end{align*}
    We show the map
    \begin{align*}
        \sigma^{\varepsilon,r}:P_{\infty}^{\varepsilon,r}(A')\rightarrow K^{2\varepsilon,r}_0(A) ;\text{ }p\mapsto [q]
    \end{align*}
    is well defined on the $K$-theoretic level,
     where $q\in M_{\infty}(A_r)$ is any element which satisfies $\|p-q\|<\frac{1}{15}\varepsilon$.
    If we have a homotopy $(p_t)_t$ throughout $P^{\varepsilon,r}_{\infty}(A')$, then there exists $k$ such that we have $\|p(t)-p(t')\|<\frac{1}{15}\varepsilon$ whenever $|t-t'|\leq \frac{1}{k}$. Take a partition of the interval
    \begin{align*}
        0=t_0<t_1<\cdots <t_{k-1}<t_k=1
    \end{align*}
    with $t_j=\frac{j}{k}$ and for each $j+0,1,\cdots, k$ we can take $q_j\in A_r$ such that $\|p_{t_j}-q_j\|< \frac{1}{20}\varepsilon$.  By linearly interpolating $q_j$'s at $t_j$, we can define another homotopy $(q_t)_t$. Then $(q_t)_t$ satisfies $\|p_t-q_t\|<\frac{1}{5}\varepsilon$ for all $t$, so especially $(q_t)_t$ is a homotopy throughout $P^{2\varepsilon,r}_{\infty}(A)$ by Remark \ref{almost same projections}. Therefore we have a map
    \begin{align*}
        \sigma^{\varepsilon,r}:K^{\varepsilon,r}_*(A')\rightarrow K^{2\varepsilon,r}_*(A);\text{ }[p,n]_{\varepsilon,r}\mapsto [q,n]_{2\varepsilon,r}.
    \end{align*}
     Clearly compositions of $\iota^{\varepsilon,r}$ and $\sigma^{\varepsilon,r}$ are the same as forgetful maps.
    }
\end{rmk}

Let 
\[
\begin{CD}
    \mathcal{O}_0 @>{\mathcal{F}}>> \mathcal{O}_1 @>{\mathcal{G}}>> \mathcal{O}_2
\end{CD}
\]
be a sequence of $(\alpha,k)$-controlled morphism $\mathcal{F}$ and $(\beta,l)$-controlled morphism $\mathcal{G}$. Assume $\mathcal{O}_i'$ is a controlled object which is controlled equivalent to $\mathcal{O}_i$ via $(\lambda,h)$-controlled equivarent morphisms
\begin{align*}
        \mathcal{H}_i&:\mathcal{O}_i \longrightarrow \mathcal{O}_i'\\
        \mathcal{J}_i&:\mathcal{O}_i' \longrightarrow \mathcal{O}_i
\end{align*} for $i=0,1,2$. We can define a $(\lambda_1,h_1)=(\lambda,h)*(\alpha,k)*(\lambda,h)$-controlled morphism $\mathcal{F}'=(F'^{\varepsilon,r})$ and $(\lambda_2,h_2)=(\lambda,h)*(\beta,l)*(\lambda,h)$-controlled morphism $\mathcal{G}'=(G'^{\varepsilon,r})$ by
\[
\begin{CD}
 F'^{\varepsilon,r}:O_0'^{\varepsilon,r} @>{J_0^{\varepsilon,r}}>> O_0^{\lambda \varepsilon, h(\varepsilon)r } @>{F^{\lambda \varepsilon, h(\varepsilon)r}}>> O_1^{\alpha \lambda \varepsilon, k(\lambda \varepsilon)h(\varepsilon)r} @>{H_1^{\alpha \lambda \varepsilon, k(\lambda \varepsilon)h(\varepsilon)r}}>> O_1'^{\alpha \lambda^2 \varepsilon,h(\alpha \lambda \varepsilon)k(\lambda \varepsilon)h(\varepsilon)r}  
\end{CD}
\]
\[
\begin{CD}
 G'^{\varepsilon,r}:O_1'^{\varepsilon,r} @>{J_1^{\varepsilon,r}}>> O_1^{\lambda \varepsilon, h(\varepsilon)r } @>{G^{\lambda \varepsilon, h(\varepsilon)r}}>> O_2^{\beta \lambda \varepsilon, l(\lambda \varepsilon)h(\varepsilon)r} @>{H_2^{\beta \lambda \varepsilon, k(\lambda \varepsilon)h(\varepsilon)r}}>> O_2'^{\beta \lambda^2 \varepsilon,h(\beta \lambda \varepsilon)l(\lambda \varepsilon)h(\varepsilon)r}  
\end{CD}
\] so that the following controlled diagram commutes;
\[
\begin{CD}
    \mathcal{O}_0 @>{\mathcal{F}}>> \mathcal{O}_1 @>{\mathcal{G}}>> \mathcal{O}_2 \\
    @A{\cong}AA  @A{\cong}AA @A{\cong}AA \\
    \mathcal{O}'_0 @>{\mathcal{F}'}>> \mathcal{O}'_1 @>{\mathcal{G}'}>> \mathcal{O}'_2 .\\  
\end{CD}
\]

We can easily see that the exactness of $(\mathcal{F},\mathcal{G})$ inherits to $(\mathcal{F}',\mathcal{G}')$.

\begin{lem}\label{equivalent-exactness}
 \rm{   If $(\mathcal{F},\mathcal{G})$ is $(\delta,p)$-exact, then there exists a controlled pair $(\delta',p')$ such that $(\mathcal{F}',\mathcal{G}')$ is $(\delta',p')$-exact.}
\end{lem}

\section{Functoriality and homotopy invariance of $K^{\varepsilon,r}_*(C^*(X))$}
In this section, we formulate how a coarse map between two locally compact metric spaces induces a map between their controlled $K$-groups and then we prove the homotopy invariance of this induced map under a certain condition. It is well known that the $K$-theory of Roe algebra is a coarse functor: if we have a coarse map between two locally compact metric spaces $X$ and $Y$, then we have an induced map
\begin{align*}
    f_*:K_*(C^*(X)) \rightarrow K_*(C^*(Y)).
\end{align*}
We give a controlled version of this morphism. After that, we discuss the homotopy invariance of this functor. This is an analogue of Lemma 4.8 in \cite{yu1998novikov}, which is a homotopy invariance of $K$-theory of a different $C^*$-algebra $C^*_{L,0}(X)$. We show the same statement for $C^*(X)$ by a similar argument.  

\begin{dfn}\label{coarse map}
\rm{
Let $(X,d_X)$ and $(Y,d_Y)$ be locally compact metric spaces and $f:X \rightarrow Y$ be a Borel map between them. For each positive number $r\geq 0$, we define the expansion function of $f$ by
\begin{align*}
    \omega_f(r):=\sup\{d_Y(f(x_1),f(x_2));d_X(x_1,x_2)<r\}.
\end{align*}
We say a map $f$ is a coarse map if $\omega_f(r)<\infty$ for each $r\geq 0$ and $f^{-1}(K)$ is a bounded set for each bounded subset $K\subset X$ of $X$.
}
    
\end{dfn}

\begin{dfn}\label{cover}
\rm{    Let $X$ and $Y$ be two locally compact metric spaces, $H_X$ and $H_Y$ are ample $X$ and $Y$-module, respectively, $f:X \rightarrow Y$ be a coarse map and $\delta>0$. An isometry 
    \begin{align*}
        V_f:H_X \longrightarrow H_Y
    \end{align*}
    is called a $\delta$-cover of $f$ if $d(y,f(x))<\delta$ for any $(y,x) \in \text{supp}(V) \subset Y \times X$.
    }
\end{dfn}

\begin{rmk}
 \rm{   For any $X,Y,f$ and $\delta$ as above, a $\delta$-cover $V_f$ of $f$ exists. This is proven in Lemma 2.4 of \cite{yu1998novikov}.}
\end{rmk}

The proof of the next theorem is provided by the Proposition 6.3.12 of \cite{higson2000analytic}.
\begin{thm}
\rm{    Under the setting of Definition \ref{cover}, for any $\delta>0$, a $\delta$-cover $V_f$ exists and $\text{Ad}(V_f)$ can restrict to 
    \begin{align*}
        \text{Ad}(V_f):C^*(X) \longrightarrow C^*(Y),
    \end{align*} and the induced map on $K$-theory 
    \begin{align*}
        f_*:=\text{Ad}(V_f)_*:K_*(C^*(X)) \longrightarrow K_*(C^*(Y))
    \end{align*}does not depend on $\delta$ and an $\delta$-cover $V_f$ of $f$.
    }
\end{thm}

\begin{rmk}
\rm{    Assume $f:X \rightarrow Y$ is a coarse map and $V_f$ is an $\delta$-cover of $f$, then we have
    \begin{align*}
     \text{prop}(\text{Ad}(V_f)(T))< \omega_f(r)+2\delta  
    \end{align*}
    for any $T \in B(H_X)$ whose propagation is finite.
    Therefore we can define 
    \begin{align*}
        \text{Ad}(V_f)_*:K^{\varepsilon,r}_*(C^*(X)) \longrightarrow K^{\varepsilon,\omega_f(r)+2\delta}_*(C^*(Y))
    \end{align*}
    }
\end{rmk}

\begin{lem}
\rm{
    Fix two positive numbers $\delta$ and $r$. Let $f:X \rightarrow Y$ be coarse maps with $\omega_f(r)<R$ and $V_f,V_f':H_X \rightarrow H_Y$ be  $\delta$-cover maps of $f$. Then we have
    \begin{align*}
        \iota^{(\varepsilon,R+2\delta),(\varepsilon,R+8\delta)}\circ \text{Ad}(V_f)_*=\iota^{(\varepsilon,R+2\delta),(\varepsilon,R+8\delta)} \circ \text{Ad}(V_f')_*,
    \end{align*}
    as maps from $K^{\varepsilon,r}_*(C^*(X))$ to $K^{\varepsilon,R+8\delta}_*(C^*(Y))$.
    }
\end{lem}

\begin{pf}
\rm{ Let $U$ be the unitary 
    \begin{align*}
        U=\begin{pmatrix}
            1-V_fV_f^* & V_fV_f'^* \\
            V_f'V_f^* & 1-V_f'V_f'^*
        \end{pmatrix} ,
    \end{align*} which is an element in the multiplier of $M_2(C^*(Y))$ and we have
\begin{align*}
    \text{Ad}U \begin{pmatrix}
        \text{Ad}(V_f)(T) & 0 \\
        0 & 0
    \end{pmatrix}
    =\begin{pmatrix}
        0 & 0 \\
        0 & \text{Ad}(V_f')(T) \\
    \end{pmatrix}
\end{align*} for any $T \in C^*(H_X)$. Since $\text{prop}(U)<2\delta$, conjugation by the unitaries 
\begin{align*}
\Tilde{U}_t:=
    \begin{pmatrix}
        U & 0 \\
        0 & 1
    \end{pmatrix}
    \begin{pmatrix}
        \cos(t) & -\sin(t) \\
        \sin(t) & \cos(t)
    \end{pmatrix}
    \begin{pmatrix}
        1 & 0 \\
        0 & U^{*}
    \end{pmatrix}
    \begin{pmatrix}
        \cos(t) & \sin(t) \\
        -\sin(t) & \cos(t)
    \end{pmatrix}\in M_4(C^*(Y)) \text{ }(t \in [0,2\pi])
\end{align*}
gives a homotopy between $\text{diag}(\text{Ad}(V_f)(p),0,0,0)$ and $\text{diag}(0,\text{Ad}(V_f')(p),0,0)$ through $(\epsilon,R+8\delta)$-projections for any $(\varepsilon,r)$-quasi-projection $p$. Therefore $[\text{Ad}(V_f)(p)]= [\text{Ad}(V_f')(p)] \in K^{\epsilon,R+8\delta}_0(C^*(Y))$. \qed
} 
\end{pf}
\begin{rmk}\label{functoriality}
\rm{  For a coarse map $f:X \rightarrow Y$, with $\omega_f(r)<R$ and any $\delta>0$ we have a well defined homomorphism:
    \begin{align*}
        \text{Ad}(V_f)_*:K^{\varepsilon,r}_*(C^*(H_X)) \longrightarrow K^{\varepsilon,R+8\delta}_*(C^*(H_Y))
    \end{align*} independently of the choice of $\delta$-cover $V_f$ of $f$.}
\end{rmk}

We use this remark in the following way. Let $H_X$ and $H_X'$ be ample $X$-modules and fix. For each $\delta>0$ there exist $\delta$-covers of $\text{id}_X$
\begin{align*}
    V_1:H_X \rightarrow H_X' \text{ and } V_2:H_X' \rightarrow H_X.
\end{align*} 
Then since $V_2V_1$ and $V_1V_2$ are $2\delta$-covers of $\text{id}_X$, by Remark \ref{functoriality}, $\text{Ad}(V_1)_*$ and $\text{Ad}(V_2)_*$ are asymptotically inverse each other in the following sense:
\begin{prop}
\rm{
Under the above setting, 
    \begin{align*}
       \text{Ad}(V_2)_* \circ \text{Ad}(V_1)_*: K^{\varepsilon,r}_*(C^*(H_X)) \rightarrow  K^{\varepsilon,r+16\delta}_*(C^*(H_X')) \rightarrow K^{\varepsilon,r+32\delta}_*(C^*(H_X))
    \end{align*}
    is equal to the forgetful map $\iota^{(\varepsilon,r),(\varepsilon,r+32\delta)}$ and the composition
    \begin{align*}
       \text{Ad}(V_1)_*\circ \text{Ad}(V_2)_*: K^{\varepsilon,r}_*(C^*(H_X')) \rightarrow  K^{\varepsilon,r+16\delta}_*(C^*(H_X)) \rightarrow K^{C\varepsilon,r+32\delta}_*(C^*(H_X'))
    \end{align*}
    is equal to the forgetful map $\iota^{(\varepsilon,r),(\varepsilon,r+32\delta)}$. 
    }
\end{prop}

Next we discuss the homotopy invariance. 
\begin{dfn}
\rm{    Two coarse maps $f,g:X \longrightarrow Y$ are said to be strongly Lipschitz homotopic if there exists a continuous homotopy
    \begin{align*}
        F:[0,1] \times X \longrightarrow Y
    \end{align*} such that }
   \begin{enumerate}[1)]
       \item $F(t,\cdot):X \longrightarrow Y$ is a proper map for each $t \in [0,1]$,

       \item the family $\{F(t,\cdot):X \longrightarrow Y\}$ is uniformly Lipschitz, namely each $F(t,\cdot)$ is Lipschitz with a Lipschitz constant $c$ that is independent on $t \in [0,1]$,

       \item $\{F(\cdot,x):[0,1] \longrightarrow Y\}_{x \in X}$ is uniformly equicontinuous, and

       \item $F(0,\cdot)=f$ and $F(1,\cdot)=g$.
       \end{enumerate}

\end{dfn}
The next theorem is almost same as Lemma 4.8 in \cite{yu1998novikov}, where the same statement for a $C^*$-algebra $C^*_{L,0}(X)$ is shown. 
\begin{thm}\label{homotopy invariance}
\rm{    Assume $f,g:X \longrightarrow Y$ be strongly Lipschitz homotopic via $F$ whose uniform Lipschitz constant is bounded by $c$ and $V_f$ and $V_g$ are their $\delta$-covers of $f$ and $g$, respectively. Then 
\begin{align*}
    \text{Ad}(V_f)_*=\text{Ad}(V_g)_*:K_*^{\varepsilon,r}(C^*(X)) \longrightarrow K_*^{21\varepsilon,5(cr+2\delta)}(C^*(Y)).
\end{align*}}
\end{thm}

\begin{pf}
\rm{   We will show this for $*=1$ and similar argument works for $*=0$. We can take a partition of the interval
    \begin{align*}
        0=t_0<t_1<\cdots<t_{l-1}<t_l=1
    \end{align*}
    such that $\|F(t_j,x)-F(t_{j+1},x)\|<\delta$ for all $x \in X$ and $j=1,2,\cdots l-1$ by the third condition. We write $f_j:=F(t_j,\cdot)$. For each $j=1,2,\cdots l$, we can take $\delta$-cover $V_j:H_X \rightarrow H_Y$ of $f_j$. Take an $(\varepsilon,r)$-quasi-projection $u$ over $\Tilde{C}^*(X)$. We define
    \begin{align*}
      u_i:=(\text{Ad}V_{f_i})(u) \text{ } (i=0,1,\cdots ,l),
    \end{align*}
    where $\text{Ad}V_{i}$ is a unital extension of $*$-homomorphism $\text{Ad}V_{i}:C^*(X) \rightarrow C^*(Y)$ to their unitizations. 
For each $i$, we define $w_i:=u_iu_l^*$, and $(3\varepsilon, 2(cr+2\delta))$-quasi-unitary operators acting on $\bigoplus_{i=0}^{l}(H_Y \oplus H_Y)$
\begin{align*}
    a:=&\bigoplus_{i=0}^{l}(w_i \oplus I)=(w_0 \oplus I)\oplus (w_1 \oplus I)\oplus \cdots \oplus (w_{l-1} \oplus I)\oplus (w_l \oplus I)\\
    b:=&\bigoplus_{i=0}^{l-1} (w_{i+1} \oplus I)\oplus (w_{l}\oplus I) =(w_1 \oplus I)\oplus (w_2 \oplus I)\oplus \cdots \oplus(w_l \oplus I) \oplus (w_l \oplus I)\\
    c:=&(w_l \oplus I)\oplus \bigoplus_{i=1}^l(w_i \oplus I)=(w_l \oplus I)\oplus(w_1 \oplus I)\oplus (w_2 \oplus I)\oplus \cdots \oplus (w_l \oplus I) .
\end{align*}
For $t \in [0,1]$, we define an isometry $V_{i,i+1}(t):(H_X \oplus H_X) \rightarrow (H_Y \oplus H_Y)$ by
\begin{align*}
    V_{i,i+1}(t)=R(t)\begin{pmatrix}
        V_i & 0 \\
        0 & V_{i+1} 
    \end{pmatrix} R(t)^*,
\end{align*} where $R(t)$ is the $\frac{\pi}{2}t$-rotation matrix $\begin{pmatrix}
    \cos(\frac{\pi}{2}t) & \sin(\frac{\pi}{2}t) \\
    -\sin(\frac{\pi}{2}t) & \cos(\frac{\pi}{2}t) \\
\end{pmatrix}$. Note that since $V_{i,i+1}(t)$ is a $2\delta$-cover of $f_i$ with respect to ample modules $(H_X\oplus H_X)$ and $(H_Y\oplus H_Y)$, $\text{Ad}(V_{i,i+1}(t))(u \oplus I)$ is $(\varepsilon,cr+4\delta)$-quasi-unitary and
\begin{align*}
    \text{Ad}(V_{i,i+1}(0))(u \oplus I) &= u_i \oplus I \\
    \text{Ad}(V_{i,i+1}(1))(u \oplus I) &= u_{i+1} \oplus I.
\end{align*}
Therefore 
\begin{align*}
    (\bigoplus_{i=0}^{l-1} \text{Ad}(V_{i,i+1}(t))(u\oplus I) \bigoplus (u_l \oplus I)) \circ \bigoplus_{i=0}^l(u_{l}^* \oplus I)
\end{align*} is a homotopy throughout $(3\varepsilon, 2(cr+4\delta))$-quasi-unitaries between $a$ and $b$. We can show that $b$ and $c$ are $(3\varepsilon, 2(cr+2\delta))$-quasi-unitary-homotopic by rotation of first $l$-coordinates and the last coordinate. So $a$ and $c$ are homotopic throughout $(3\varepsilon, 2(cr+4\delta))$-quasi-unitaries. 

Note that since $\|u_0-w_l^*w_0u_l\|\leq 4\varepsilon$, two $(21\varepsilon,5(cr+4\delta))$-quasi-unitaries 
\begin{align*}
     (u_0\oplus I) \oplus \bigoplus_{k=1}^l(I \oplus I) \text{ and }c^*a((u_{l}\oplus I) \oplus \bigoplus_{k=1}^l(I \oplus I))
\end{align*}
are homotopic throughout $(21\varepsilon,5(cr+4\delta))$-quasi-unitaries by Remark \ref{almost same projections}. Similarly, 
\begin{align*}
    a^*a((u_{l}\oplus I) \oplus \bigoplus_{k=1}^l(I \oplus I)) \text{ and }
    (u_{l}\oplus I) \oplus \bigoplus_{k=1}^l(I \oplus I)
\end{align*}
are homotopic throughout $(21\varepsilon,5(cr+4\delta))$-quasi-unitaries.
We can construct a homotopy throughout $(6\varepsilon, 5(cr+4\delta))$-quasi-unitaries between 
\begin{align*}
    c^*a((u_{l}\oplus I) \oplus \bigoplus_{k=1}^l(I \oplus I))\text{ and }
    a^*a((u_{l}\oplus I) \oplus \bigoplus_{k=1}^l(I \oplus I))
\end{align*}
by using the homotopy between $a$ and $c$.
This proves $$\text{Ad}(V_f)_*([u])=\text{Ad}(V_g)_*([u]) \in K_1^{21\varepsilon,5(cr+4\delta)}(C^*(H_Y)).$$
\qed
}
\end{pf}

\section{Asymptotic Mayer-Vietoris exact sequence of controlled $K$-theory}
In this section, we decompose an $n$-dimensional simplicial complex $X$ into two pieces to obtain a Mayer-Vietoris sequence of controlled $K$-theories. First, we recall the condition of decomposition of filtered $C^*$-algebra to obtain a controlled Mayer-Vietoris sequence following \cite{oyono2019quantitative} and then we show that this decomposition of a simplicial complex satisfies the conditions. Our purpose in this section is to build an asymptotically long exact sequence to obtain information of $K^{\varepsilon,r}_*(C^*(X))$ from finite simplicial complexes whose dimension is less than $\text{dim}(X)$. We encode the spherical metric on each simplex (Definition 7.2.1 of \cite{willett2020higher}).

To obtain a Mayer-Vietoris sequence in terms of two subspaces (or subalgebras), these two have to intersect enough. We state these conditions.

\begin{dfn}[\cite{oyono2019quantitative} Definition 2.6 and Remark 2.7] \label{completely coercive decomposition}
\rm{
    Let $A=(A_s)_s$ be a filtered $C^*$-algebra and let $r$ be a positive number. A completely coercive decomposition pair of degree $r$ for $A$ is a pair $(\Delta_1,\Delta_2)$ of closed linear subspaces of $A_r$ such that there exists a positive number $C$ satisfying the following:

    for any positive number $s$ with $s \leq r$, any positive integer $n$ and any $x \in M_n(A_s)$, there exist $y\in M_n(A_s\cap \Delta_1)$ and $z \in M_n(A_s\cap\Delta_2)$ with
    \begin{align*}
        \|y\|\leq C\|x\| \text{, }\|z\|\leq C\|x\| \text{ and }x=y+z.
    \end{align*}
    The constant $C$ is called the coercity of the decomposition.
    }
\end{dfn}

\begin{dfn}[\cite{oyono2019quantitative} Definition 2.8 and Definition 2.13]
\rm{
    Let $A=(A_s)_s$ be a filtered $C^*$-algebra, $r$ be any positive numbers and $\Delta$ be a closed subspace of $A_r$. A filtered $C^*$-subalgebra $B=(B\cap A_s)_s$ of $A$ is called an $r$-neighborhhod of $\Delta$ if it contains the subspace
    \begin{align*}
        N_{\Delta}^{(r,5r)}=\Delta+A_{5r}\cdot \Delta+\Delta\cdot A_{5r}+ A_{5r}\cdot\Delta\cdot A_{5r}.
    \end{align*}
    }
\end{dfn}

\begin{dfn}[\cite{oyono2019quantitative} Definition 2.15]\label{CIA}
\rm{
Let $A$ be a $C^*$-algebra. A pair $(S_1,S_2)$ of subsets of $A$ is said to satisfy complete intersection approximation property (CIA) if there exists $C>0$ such that for any positive number $\varepsilon$, any positive integer $n$ and any $x\in M_n(S_1)$ and $M_n(S_2)$ with $\|x-y\|\leq\varepsilon$, there exists $z \in M_n(S_1 \cap S_2)$ such that
\begin{align*}
    \|z-x\|\leq C\varepsilon\text{, } \|z-y\|\leq C\varepsilon.
\end{align*}
The constant $C$ is called the coercity of the pair $(S_1,S_2)$.
    }
\end{dfn}

\begin{dfn}[\cite{oyono2019quantitative} Definition 2.16]\label{weak Mayer-Vietoris pair}
    \rm{
    Let $A=(A_s)_s$ be a filtered $C^*$-algebra and $r$ be any positive number. An $r$-controlled weak Mayer-Vietoris pair for $A$ is a quadruple $(\Delta_1,\Delta_2,A_{\Delta_1},A_{\Delta_2})$ such that for a positive number $C>0$,
    \begin{enumerate}[1)]
        \item $(\Delta_1,\Delta_2)$ is a completely coercive decomposition pair for $A$ of order $r$ with coercity $C$,
        \item $A_{\Delta_i}$ is an $r$-controlled neighborhood of $\Delta_i$ for $i=1,2$,
        \item the pair $(A_{\Delta_1,s},A_{\Delta_2,s})$ has the CIA property for any positive number $s\leq r$ with coercity $C$.
    \end{enumerate}
    The number $C$ is called the coercity of the $r$-controlled Mayer-Vietoris pair $(\Delta_1,\Delta_2,A_{\Delta_1},A_{\Delta_2})$.
    }
\end{dfn}

Now we can state the main technical tool of this section.

\begin{thm}[\cite{oyono2019quantitative} Theorem 3.10]\label{Mayer-Vietoris}
\rm{
For any positive number $C$, there exists a control pair $(\lambda,h)$ such that for any filtered $C^*$-algebra $A=(A_s)_s$, any positive number $r$ and any $r$-controlled weak Mayer-Vietoris pair $(\Delta_1,\Delta_2,A_{\Delta_1},A_{\Delta_2})$ for $A$ of order $r$ with coercity $C$, we have a $(\lambda,h)$-exact six term exact sequence of order $r$:

\[
\begin{CD}
    \mathcal{K}_0(A_{\Delta_1} \cap A_{\Delta_2}) @>>> \mathcal{K}_0(A_{\Delta_1}) \bigoplus \mathcal{K}_0(A_{\Delta_2}) @>>> \mathcal{K}_0(A) \\
    @A{\partial}AA @. @V{\partial}VV \\
    \mathcal{K}_1(A_{\Delta_1} \cap A_{\Delta_2}) @>>> \mathcal{K}_1(A_{\Delta_1}) \bigoplus \mathcal{K}_1(A_{\Delta_2}) @>>> \mathcal{K}_1(A) .
\end{CD}
\]

}
\end{thm}

Let $X$ be a $n$-dimensional finite simplicial complex. For each $n$-simplex $Y$ in $X$, we define 
\begin{align*}
    Y_1:=\{x \in Y;d(c_Y,x)\leq\frac{1+\frac{1}{10}}{2}\}
\end{align*}
\begin{align*}
    Y_2:=\{x \in Y;d(c_Y,x)\geq\frac{1-\frac{1}{10}}{2}\},
\end{align*} where $c_Y$ is the center of each simplex. We decompose $X$ into two subsets
\begin{align*}
    X_1 := \bigcup \{Y_1;\text{Y is an $n$-dimensional simplex in $X$}\}  
\end{align*}
and
\begin{align*}
    X_2 := \bigcup \{Y_2;\text{Y is an $n$-dimensional simplex in $X$}\} \cup \bigcup\{Y \in X; \dim Y<n \}.
\end{align*}
Then $(X_1)_{\frac{1}{10}}$, $(X_2)_{\frac{1}{10}}$ and $(X_1 \cap X_2)_{\frac{1}{10}}$ are strongly Lipschitz homotopic to a finite $0$-dimensional complex $Z_1$, $(n-1)$-dimensional complex $Z_2$ and $Z$, respectively with a uniform Lipschitz constant $c_n$ depending only on $n$, where $(W)_r$ is an $r$-neighborhood of $W \subset X$ in $X$. We denote these strong Lipschitz homotopy equivelence maps by;
\begin{align*}
    f_1 : (X_1)_{\frac{1}{10}} \longrightarrow Z_1, g_1 : Z_1 \longrightarrow (X_1)_{\frac{1}{10}}
\end{align*}
\begin{align*}
    f_2 : (X_2)_{\frac{1}{10}} \longrightarrow Z_2, g_2 : Z_2 \longrightarrow (X_2)_{\frac{1}{10}}
\end{align*}
\begin{align*}
    f : (X_1\cap X_2)_{\frac{1}{10}} \longrightarrow Z, g:Z \longrightarrow (X_1\cap X_2)_{\frac{1}{10}}.
\end{align*} 
By Theorem \ref{homotopy invariance}, $(f_1,g_1)$, $(f_2,g_2)$ and $(f,g)$ induce asymptotically inverse maps each other on controlled $K$-theories.
We denote
\begin{align*}
    \Delta_i:=C^*(X_i), A_i:=C^*((X_i)_{\frac{1}{10}})
\end{align*} for each $i=1,2$. We show the quadruple $(\Delta_1,\Delta_2,A_1,A_2)$ satisfies the condition to be a $\frac{1}{50}$-controlled Mayer-Vietoris pair for $C^*(X)$. 

(1) Show the pair $(\Delta_1,\Delta_2)$ is a completely coercive decomposition pair of defree $\frac{1}{50}$. Take any positive number $r$ with $0<r<\frac{1}{50}$ and $x \in (C^*(X))_r$. We denote $\Pi_1=X_1 \setminus (X_1 \cap X_2)$, $\Pi_2=X_1 \cap X_2$ and $\Pi_3=X_2 \setminus (X_1 \cap X_2)$. In terms of this disjoint decomposition we write $x=(x_{i,j})_{1\leq i,j \leq 3}$, where $x_{i,j}=\chi_{\Pi_i}x\chi_{\Pi_j}$. Since $x$ can be approximated by operators whose propagation is smaller than $r$, $x$ has the form
\begin{align*}
    x=\begin{pmatrix}
        x_{1,1} & x_{1,2} & 0 \\
        x_{2,1} & x_{2,2} & x_{2,3} \\
        0 & x_{3,2} & x_{3,3}
    \end{pmatrix}
\end{align*}
we define 
\begin{align*}
    x_1=\begin{pmatrix}
        x_{1,1} & x_{1,2} & 0 \\
        x_{2,1} & x_{2,2} & 0 \\
        0 & 0 & 0
    \end{pmatrix}
    \text{, }
    x_2=\begin{pmatrix}
        0 & 0 & 0 \\
        0 & 0 & x_{2,3} \\
        0 & x_{3,2} & x_{3,3}
    \end{pmatrix}.
\end{align*} Then we have $x_i \in \Delta_i$, $x=x_1+x_2$, $\|x_i\|\leq 4 \|x\|$ for $i=1,2$. If $x \in M_n((C^*(X))_r)$, we can regard $x$ is an operator on $\oplus_{k=1}^nH_X$ with the diagonal module structure so that we can apply the same argument.

(2) Clearly $A_i$ is an ${\frac{1}{50}}$-controlled $\Delta_i$-neighborhood-$C^*$-algebra. 

(3) Show the pair $((A_1)_r,(A_2)_r)$ has the CIA property for any $0 \leq r \leq \frac{1}{50}$. Let $\varepsilon >0, x \in (A_1)_r$ and $y \in (A_2)_r$ with $\|x-y\|<\varepsilon$. We denote $\Sigma_1=(X_1)_{\frac{1}{10}+r} \setminus ((X_1)_{\frac{1}{10}+r} \cap (X_2)_{\frac{1}{10}+r})$, $\Sigma_2=(X_1)_{\frac{1}{10}+r} \cap (X_2)_{\frac{1}{10}+r}$ and $\Sigma_3=(X_2)_{\frac{1}{10}+r} \setminus ((X_1)_{\frac{1}{10}+r} \cap (X_2)_{\frac{1}{10}+r})$. In terms of this disjoint decomposition we write $x=(x_{i,j})_{1\leq i,j \leq 3}$ and $y=(y_{i,j})_{1\leq i,j \leq 3}$, where $x_{i,j}=\chi_{\Sigma_i}x\chi_{\Sigma_j}$ and similarly for $y$. Then,
\begin{align*}x-y=
    \begin{pmatrix}
        x_{1,1} & x_{1,2} & 0 \\
        x_{2,1} & x_{2,2}-y_{2,2} & -y_{2,3} \\
        0 & -y_{3,2} & -y_{3,3}
    \end{pmatrix}. 
\end{align*}
Define 
\begin{align*}
  z=\frac{x_{2,2}+y_{2,2}}{2} = 
  \begin{pmatrix}
  0&0&0\\
  0 & \frac{x_{2,2}+y_{2,2}}{2} & 0 \\
  0&0&0\\
  \end{pmatrix}\in (A_1)_r \cap (A_2)_r.  
\end{align*}
Then 
\begin{align*}
    \|x-z\|=
   \left\| \begin{pmatrix}
       x_{1,1} & x_{1,2} & 0 \\
       x_{2,1} & 0 & 0 \\
       0 & 0 & 0 
    \end{pmatrix}+
    \begin{pmatrix}
  0&0&0\\
  0 & \frac{x_{2,2}-y_{2,2}}{2} & 0 \\
  0&0&0\\
  \end{pmatrix}\right\| \leq 4\varepsilon
\end{align*} and similarly for $y$. If $x$ and $y$ are matrices, we can reduce it to the case we have just proven as we did in (1). Therefore the pair $((A_1)_\frac{1}{10},(A_2)_\frac{1}{10})$ has CIA property with coercity $4$. Therefore by Theorem 3.10 in \cite{oyono2019quantitative}, there exists a controll pair $(\delta,p)$ such that we have the following $(\delta,p)$ exact sequence of order $\frac{1}{50}$:
\[
\begin{CD}
    \mathcal{K}_1(A_1 \cap A_2) @>{k}>> \mathcal{K}_1(A_1) \bigoplus \mathcal{K}_0(A_2) @>{l}>> \mathcal{K}_1(C^*(X)) \\
    @A{\partial}AA @. @V{\partial}VV \\
    \mathcal{K}_0(C^*(X)) @<{l}<< \mathcal{K}_0(A_1) \bigoplus \mathcal{K}_0(A_2) @<{k}<< \mathcal{K}_0(A_1 \cap A_2)
\end{CD}
\]

By Lemma \ref{equivalent-exactness} and lemma \ref{homotopy invariance}, there exist control pairs $(\lambda_n,h_n)$ and $(\delta'_n,p'_n)$ depending only on $n$ such that their exist $(\lambda_n, h_n)$-morphisms $k$, $l$ and $\partial$ in  the following sequence, which are $(\delta'_n,p'_n)$-exact of order $s_n$ depending only on $n$:

\[
\begin{CD}
    \KK_1(C^*(Z)) @>{k}>> \KK_1(C^*(Z_1)) \bigoplus \KK_0(C^*(Z_2)) @>{l}>> \KK_1(C^*(X)) \\
    @A{\partial}AA @. @V{\partial}VV \\
    \KK_0(C^*(X)) @<{l}<< \KK_0(C^*(Z_1)) \bigoplus \KK_0(C^*(Z_2)) @<{k}<< \KK_0(C^*(Z))
\end{CD}
\]

\section{Quantitative description of $K$-homology}

In this section, we show the main step of the main theorem, which states that the $K$-homology $K_*(X)$ of a finite simplicial complex $X$ can be realized as an image of a forgetful map of the controlled $K$-theory $K^{\epsilon,r}_*(C^*(X))$. We show this using the asymptotic version of the five lemma between Mayer-Vietoris sequences of $K$-homology and controlled $K$-theory. For this purpose, it is convenient to reformulate the Mayer-Vietoris sequence of $K$-homology in terms of controlled $K$-theory so that the diagram commutes. For the next lemma, we give the definition of the localization algebra $C^*_L(X)$, which is defined in \cite{yu1997localization}.

\begin{dfn}\label{localization}
\rm{
For an ample $X$-module $H_X$, we define the algebraic localization algebra $C^*_{L,\text{alg}}(H_X)$ to be the algebra consists of uniformly continuous and bounded functions from $[1,\infty)$ to the algebraic Roe algebra $C^*_{alg}(H_X)$ and the localization algebra $C^*_L(H_X)$ to be the completion of $C^*_{L,\text{alg}}(X)$ by the sup-norm. Again, when it is clear from the context, we will write $C^*_{L,\text{alg}}$ and $C^*_L(X)$ in place of $C^*_{L,\text{alg}}(H_X)$and $C^*_L(H_X)$, respectively. Localization algebra is also filtered by the supremum of propagations i.e.
\begin{align*}
    C^*_{L}(X)_r=\{f\in C^*_{L,\text{alg}}(X):\sup_{t\in [1,\infty)} \text{prop}(f(t))<r\}
\end{align*}
or its closure. (See Remark \ref{closed or non-closed}.)
  }   
\end{dfn}

Then Yu showed that the $K$-theory of the localization algebra is isomorphic to $K$-theory for any finite dimensional simplicial complex.

\begin{thm}[\cite{yu1997localization} Theorem 3.2]
\rm{
    If $X$ is a finite dimensional simplicial complex endowed with a spherical metric, then we have an isomorphism $K_*(X)\cong K_*(C^*_L(X))$.
    }
\end{thm}
To use the same Mayer-Vietoris sequence for $K$-homology and controlled $K$-theory, we formulate $K$-homology in terms of controlled $K$-theory.

\begin{lem}\label{k-homology}
\rm{
    For any $0<\varepsilon<\frac{1}{8}$ and $r>0$, we have
    \begin{align*}
        K_*(C^*_L(X))=K_*^{\varepsilon,r}(C^*_L(X)).
    \end{align*}
 }   
\end{lem}

\begin{pf}
    \rm{
Since $C^*_{L,\text{alg}}(X)$ is dense in $C^*_L(X)$, by Proposition 4.9 of \cite{guentner2016dynamical}, we have an isomorphism  
\begin{align*}
     K_*^{\varepsilon}(C^*_{L,\text{alg}}(X)) \cong K^{\epsilon}_*(C^*_L(X))\cong K_*(C^*_L(X)) .
\end{align*}
    So it suffices to construct an isomorphism
\begin{align*}
    K_*^{\varepsilon,r}(C^*_L(X)) \cong K_*^{\varepsilon}(C^*_{L,\text{alg}}(X)).
\end{align*}    
We construct an inverse of the forgetful map $\iota^{\varepsilon,r}: K_*^{\varepsilon,r}(C^*_L(X)) \rightarrow K_*^{\varepsilon}(C^*_{L,\text{alg}}(X))$. Take any $\varepsilon$-projection $(p_t)_{t\in [1,\infty)}$ over $C^*_{L,\text{alg}}(X)$. By the definition of $C^*_{L,\text{alg}}(X)$, there exists $N \in [1,\infty)$ such that 
\begin{align}\label{large N}
    \text{prop}(p_t)<r \text{ for all } t \geq N .
\end{align}
We define
\begin{align*}
    \sigma^{\varepsilon,r}:P^{\epsilon,r}_{\infty}(C^*_{L,\text{alg}}(X)) \rightarrow K_*^{2\varepsilon,r}(C^*_L(X)) \text{ ; } [(p_t)_t] \mapsto [(p_{N+t})_t]
\end{align*}
independently of the choice of $N$ satisfying \eqref{large N}.
We show that this is well defined on the $K$-theoretic level. Assume we have a homotopy $\{(p(s)_t)_t\}_{s \in [0,1]}$ in $\varepsilon$-projections over $C^*_{L,\text{alg}}(X)$ parametrized by $s$. We show $\sigma^{\varepsilon,r}([p(0)_t])=\sigma^{\varepsilon,r}([p(1)_t])$. The issue is that we can not take $N$ uniformly for $s$ in general. We take $k$ such that $\|(p(s)_t)_t-(p(s')_t)_t\|<\frac{1}{15}\varepsilon$ whenever $\|s-s'\|<\frac{1}{k}$ and a partition of the interval
\begin{align*}
    0=s_0<s_1<\cdots<s_{k-1}<s_k=1,
\end{align*}
with $s_j=\frac{j}{k}$ for $j=0,1,\cdots k$. For each $j$, we can take $N_j$ such that
\begin{align*}
    \text{prop}(p(s_j)_t)<r \text{ for all } t \geq N_j .
\end{align*} We define a new path $q(s)_t$ by linearly connecting between $p(s_j)$'s. Then by Remark \ref{almost same projections}, $q(s)$ is a path through $2\varepsilon$-projections and 
\begin{align*}
    \text{prop}(q(s)_t)<r \text{ for all } t \geq M:=\max\{N_j\}_{j=1}^k .
\end{align*}
Therefore  $\sigma^{\varepsilon,r}([p(0)_t])=\sigma^{\varepsilon,r}([p(1)_t]) \in  K_*^{2\varepsilon,r}(C^*_L(X))$. So we have a map on $K$-theory, which is still denoted by
\begin{align*}
    \sigma^{\varepsilon,r}:K_*^{\epsilon}(C^*_{L,\text{alg}}(X)) \rightarrow K_*^{2\varepsilon,r}(C^*_L(X)) \text{ ; } [(p_t)_t] \mapsto [(p_{N+t})_t].
\end{align*}
The statement follows from the following commutative diagram:
$$
\begin{tikzcd}
 K_*^{\varepsilon}(C^*_{L,\text{alg}}(X))   \\ 
 K_*^{\frac{\varepsilon}{2}}(C^*_{L,\text{alg}}(X)) \arrow[u, "\cong"] \arrow[r,"\sigma^{\frac{\epsilon}{2},r}"] & K_*^{\varepsilon,r}(C^*_L(X)) \arrow[lu, "\iota^{\epsilon,r}" pos=0.5, swap].
\end{tikzcd}
$$
\qed
}
\end{pf}

\begin{dfn}\label{relaxed}
\rm{    Let $X$ be a locally compact metric space, $(\alpha,k)$ be a control pair. We define a new quantitative object $\mathcal{K}_*^{(\alpha,k)}(C^*(X))=(K_*^{(\alpha,k),(\varepsilon,r)}(C^*(X)))_{\varepsilon,r}$ by
    \begin{align*}
        K_*^{(\alpha,k),(\varepsilon,r)}(C^*(X)) = \iota ^{\varepsilon, \alpha \varepsilon, r, k(\varepsilon)r}(K_*^{\varepsilon,r}(C^*(X))) \subset K_*^{\varepsilon\alpha,k(\varepsilon)r}(C^*(X)).
    \end{align*}
    We call the group $K_*^{(\alpha,k),(\varepsilon,r)}(C^*(X))$ the $(\alpha,k)$-relaxed $(\varepsilon,r)$-controlloed $K$-theory of $C^*(X)$ because it is generated by the same generator as $K^{\varepsilon,r}_*(C^*(X))$ with a relaxed homotopy relation.
    Let $Y$ be another locally compact metric space and $\mathcal{F}=(F^{\varepsilon,r})$ be a $(\lambda,h)$-controlled morphism from $(K^{\varepsilon,r}_*(C^*(X)))_{\varepsilon,r}$ to $(K^{\varepsilon,r}_*(C^*(Y)))_{\varepsilon,r}$. Then we define a control pair $(\lambda,h')$ by
    \begin{align*}
        h'(\varepsilon)=h(\alpha\varepsilon)k(\varepsilon)/k(\lambda \varepsilon)
    \end{align*} and a $(\lambda,h')$-control morphism $\mathcal{F}'=(F'^{\varepsilon,r})$
    to be the restriction of
    \begin{align*}
        F^{\alpha \varepsilon,k(\varepsilon)r}: K_*^{\alpha \varepsilon,k(\varepsilon)r}(C^*(X)) \rightarrow K_*^{\lambda \alpha \varepsilon, h(\alpha \varepsilon)k(\varepsilon)r }(C^*(Y)) 
    \end{align*}}    from $K_*^{(\alpha,k),(\varepsilon,r)}(C^*(X))$ to $K_*^{(\alpha,k),(\lambda \varepsilon,h'(\varepsilon)r)}(C^*(Y))$.
\end{dfn}

\begin{rmk}\label{relaxed exactness}
\rm{    The $(\delta,l)$ exactness of $\mathcal{K}_*(\cdot)$ inherits to the $(\delta, l')$ exactness of $\mathcal{K}_*^{(\alpha,k)}(\cdot)$, where
    \begin{align*}
        l'(\varepsilon)=l(\alpha \varepsilon) k(\varepsilon)/k(\delta \varepsilon).
    \end{align*}}
\end{rmk}

Here we can prove the main theorem: the quantitative description of $K$-homology.

\begin{thm}
 \rm{   For each $n \in \mathbb{N}$, there exists a control pair $(\lambda_n, h_n)$ and a pair of positive numbers $(\epsilon_n,r_n)$ depending only on $n$ such that for any $n$ dimensional finite simplicial complex $X$ endowed with a spherical metric, we have 
    \begin{align*}
        K_*(X) \cong K_*^{(\lambda_n,h_n),(\varepsilon,r)}(X)
    \end{align*}
    for any $(\varepsilon,r)<(\varepsilon_n,r_n)$.}
\end{thm}

\begin{pf}
\rm{    We show this by induction on $n$. First assume $n=0$, in this case $X$ is a set of finitely many points $X=\{p_1,p_2,...p_N\}$. Let $r_0$ be the minimum of the distances of each two points;
    \begin{align*}
        r_0:= \min\{d(p_i,p_j);1\leq i,j\leq N, i \neq j\}.
    \end{align*}
    (Note that in the spherical metric, the distance between two points in different connected components is infinity but here we care the case that the distance between them can be finite because an intersection of two subcomplex can have such metric. But the distances between such two points are bounded below by $1$.)
    Since any operator whose propagation is smaller than $r_0$ have to be supported on a single point, for any $0<\varepsilon<\frac{1}{4}$ and $0<r < r_0$, we have
\begin{align*}
    K^{\varepsilon,r}_*(C^*(X))= \bigoplus_j K^{\varepsilon,r}_*(C^*(\{p_j\})) = \bigoplus_j K_*(C^*(\{p_j\})) = \bigoplus_j K_*(\mathcal{K}(H)) = \bigoplus_j K_*(\CC),
\end{align*} which is the same as $K$-homology group of $X=\{p_1,p_2,...,p_n\}$ by the cluster axiom of $K$-homology. In this case we can take $(\lambda,h)=(1,1)$ to show 
\begin{align*}
         K_*^{(\lambda_0,h_0),(\varepsilon,r)}(X) \cong K_*(X)
\end{align*} for $0<r<r_0$.
Next we assume that the statement holds for $\dim X=1,2,...,n-1$. By the previous section and Remark \ref{relaxed exactness}, we have the following diagram of $(\delta_n,l_n)$-exact sequence of $(\alpha_n,k_n)$-morphisms for some control pair $(\delta_n,l_n)$ and $(\alpha_n,k_n)$ of degree $s_{n-1}$

\[
\begin{CD}
   \cdots \rightarrow K_*(Z) @>>> K_*(Z_1) \bigoplus K_*(Z_2) \rightarrow \\
   @VVV @VVV \\
   \cdots \rightarrow \KK^{(\lambda_{n-1},h_{n-1})}_*(Z) @>>> \KK^{(\lambda_{n-1},h_{n-1})}_*(Z_1) \bigoplus \KK^{(\lambda_{n-1},h_{n-1})}_*(Z_2) \rightarrow \\
   \rightarrow K_*(X) @>>> K_{*+1}(Z) \rightarrow \cdots\\
   @V{\partial^{n-1}}VV @VVV \\
   \rightarrow \KK^{(\lambda_{n-1},h_{n-1})}_*(X) @>>> \KK^{(\lambda_{n-1},h_{n-1})}_{*+1}(Z) \rightarrow \cdots. \\
\end{CD}
\]
Note that the Mayer-Vietoris sequence of $K$-homology comes from that of controlled $K$-theory via the isomorphism in Lemma \ref{k-homology}.
The vertical maps are induced by the evaluation maps at $1$
\begin{align*}
    \partial_{(\varepsilon,r)}^{n-1}:K_*(Y)=K_*^{\varepsilon,r}(C^*_L(Y)) \rightarrow K_*^{(\lambda_{n-1},h_{n-1}),(\varepsilon,r)}(Y)
\end{align*}
for $Y=X,Z_1,Z_2,Z$.
Therefore the above diagram is commutative because of the construction after Remark 3.4 of \cite{oyono2019quantitative}.
By induction hypothesis there exists $(\varepsilon_{n-1},r_{n-1})$ such that $K^{(\lambda_{n-1},h_{n-1}),(\varepsilon,r)}_*(Y)$ is stable and isomorphic to $K_*(Y)$ if $(\varepsilon,r)<(\varepsilon_{n-1},r_{n-1})$ and $r<s_{n-1}$ for $Y=Z,Z_1,Z_2$. We take $(\varepsilon,r)$ such that 
\begin{align}\label{assumption 1}
    (\lambda_{n-1},h_{n-1})(\varepsilon,r)<(\varepsilon_{n-1},r_{n-1})
\end{align} and
\begin{align}\label{assumption 2}
    (\alpha_n,k_n)(\varepsilon,r)<(\varepsilon_{n-1},r_{n-1}).
\end{align}
 Then by diagram chasing we can show that for any $y \in K^{(\lambda_{n-1},h_{n-1}),(\varepsilon,r)}_*(X)$, there exists $x \in K_*(X)$
such that 
\begin{align*}
    \partial_{(\alpha_n \delta_n \varepsilon, K_n(\delta_n \varepsilon)l_n(\varepsilon)r)}^{n-1}(x)=\iota^{(\varepsilon,r),(\alpha_n,k_n)*(\delta_n,l_n)(\varepsilon,r)}(y).
\end{align*}
So if we define $(\lambda_n,h_n):=(\alpha_n,k_n)*(\delta_n,l_n)*(\lambda_{n-1},h_{n-1})$, 
\begin{align*}
   \partial_{(\varepsilon,r)}^{n}:K_*(X) \rightarrow K^{(\lambda_n,h_n),(\varepsilon,r)}_*(X)
\end{align*}
is surjective for all $(\varepsilon,r)$ with \eqref{assumption 1} and \eqref{assumption 2}. 
We show $\partial_{(\varepsilon,r)}^{n}$ is injective for small $(\varepsilon,r)$.
Let $(\varepsilon',r'):=(\frac{\lambda_n}{\lambda_{n-1}}\varepsilon, \frac{h_n(\varepsilon)}{h_{n-1}(\frac{\lambda_n}{\lambda_{n-1}}\varepsilon)}r)$ so that it satisfies
\begin{align}\label{from n-1 to n}
    (\lambda_{n-1} \varepsilon',h_{n-1}(\varepsilon')r')=(\lambda_n \varepsilon, h_n(\varepsilon)r).
\end{align}
Note that we have a commutative diagram
\[
\begin{CD}
    K_*(Y) @>{\partial_{\varepsilon,r}}>> K_*^{\varepsilon,r}(C^*(Y)) \\
    @V{\partial_{\varepsilon',r'}}VV @V{\iota^{\varepsilon,\lambda_n \varepsilon, r, h_n(\varepsilon)r}}VV \\
    K_*^{\varepsilon',r'}(C^*(Y)) @>{\iota^{\varepsilon', \varepsilon' \lambda_{n-1}, r', h_{n-1}(\varepsilon_1)r'}}>> K_*^{\lambda_n \varepsilon,h_n(\varepsilon)r}(C^*(Y))    .
\end{CD}
\]
 
By the induction hypothesis if we choose $(\varepsilon',r')$ small enough, then for $Y= Z,Z_1,Z_2$, the composition $\iota^{\varepsilon', \varepsilon' \lambda_{n-1}, r', h_{n-1}(\varepsilon)r'} \circ \partial_{\varepsilon',r'}$ is an isomorphism into $K_*^{(\lambda_n,h_n),(\varepsilon,r)}(Y) \subset K_*^{\lambda_n \varepsilon,h_n(\varepsilon)r}(C^*(Y)) $ by \eqref{from n-1 to n}. Therefore 
\begin{align*}
    \partial_{(\varepsilon,r)}^n : K_*(Y) \rightarrow K_*^{(\lambda_n,h_n),(\varepsilon,r)}(Y)
\end{align*}
is also isomorphic for $(\epsilon,r)$ with
\begin{align}\label{assumption 3}
    (\epsilon',r')<(\varepsilon_{n-1},r_{n-1}).
\end{align} 
 besides \eqref{assumption 1} and \eqref{assumption 2}.
Then by the following diagram chasing of commutative asymptotic exact sequence of degree $s_n$
\[
\begin{CD}
   \cdots \rightarrow K_*(Z) @>>>  K_*(Z_1) \bigoplus K_*(Z_2) \rightarrow \\
   @VVV @VVV \\
   \cdots \rightarrow \KK^{(\lambda_n,h_n)}_*(Z) @>>> \KK^{(\lambda_n,h_n)}_*(Z_1) \bigoplus \KK^{(\lambda_n,h_n)}_*(Z_2) \rightarrow \\
   \rightarrow K_{*+1}(X) @>>> K_{*+1}(Z) \rightarrow \cdots\\
   @VVV @VVV \\
   \rightarrow \KK^{(\lambda_n,h_n)}_{*+1}(X) @>>> \KK^{(\lambda_n,h_n)}_{*+1}(Z) \rightarrow \cdots,\\ 
\end{CD}
\]
we can show 
\begin{align*}
    \partial : K_*(X) \rightarrow K_{*}^{(\lambda_n,h_n),(\varepsilon,r)}(X)
\end{align*}
is also injective for sufficiently small $(\varepsilon,r)$ depending only on $n$.
\qed

\begin{rmk}\label{discretization}
\rm{
We continue assuming $X$ is a finite $n$-dimensional simplicial complex.
Let $Y$ be a countable dense subset of $X$ which is an inductive limit of finite subsets $(Y_i)$ of $Y$ and $H$ be any infinite dimensional separable Hilbert space.  We take an increasing sequence $(H_i)_i$ of finite dimensional subspaces $H_i$ of $H$ such that $H=\overline{\cup_i H_i}$ and we denote by $Q_i \in B(\ell^2(Y)\otimes H)$ the orthogonal projection onto $\ell^2(Y_i)\otimes H_i$.
We have a unital ample representation of $C(X)$ to $\ell^2(Y)\otimes H$ by a multiplication to $\ell^2(Y)$ and the identity for $H$ and this representation can be restricted to each $\ell^2(Y_i)\otimes H_i$.

In this case, the Roe algebra $C^*_{\rho}(\ell^2(Y)\otimes H)$ is the set of compact operators $\mathcal{K}(\ell^2(Y)\otimes H)=\varinjlim B(\ell^2(Y_i)\otimes H_i)$. Therefore for any $(\varepsilon,r)$-projection $p$ over $C^*_{\rho}(\ell^2(Y)\otimes H)$ and $\delta>0$, we can find $i$ such that $\|p-Q_ipQ_i\|\leq \delta$. By Remark \ref{almost same projections}, $p$ and $Q_ipQ_i$ give the same element in $K_{0}^{(\lambda_n,h_n),(\varepsilon,r)}(X)$ for sufficiently small $\delta$. 
}
\end{rmk}

}
\end{pf}

\section*{Acknowledgement}
The author would like to thank his advisor Prof. Guoliang Yu for the discussions and encouragements.

\bibliographystyle{plain}
\bibliography{main}

\end{document}